\documentstyle[11pt]{article}

\newtheorem{theorem}{Theorem}
\newtheorem{lemma}[theorem]{Lemma}
\newtheorem{proposition}[theorem]{Proposition}
\newtheorem{claim}[theorem]{Claim}
\newtheorem{conjecture}[theorem]{Conjecture}
\newtheorem{corollary}[theorem]{Corollary}
\newtheorem{definition}{Definition}

\def\qedf{\hfill $\Box$}

\large

\title{{\bf  Approximate Multipartite Version of the Hajnal--Szemer\'edi Theorem}}

\author{
B\'ela Csaba\thanks{Part of this work was done while the author worked at the Analysis and Stochastics
Research Group at the University of Szeged. Partially supported by OTKA T049398. e-mail: bela.csaba@wku.edu} \\
Department of Mathematics\\ Western Kentucky University\\
\and Marcelo Mydlarz\thanks{e-mail: marcem@cs.rutgers.edu} \\ Department of Computer Science \\ Rutgers University}

\date{}
\begin{document}
\maketitle

\begin{abstract}
Let $q$ be a positve integer, and $G$  be a $q$-partite simple graph on $qn$ vertices, with $n$ vertices in each 
vertex class. Let $\delta={k_q \over k_q+1}$,  where $k_q=q+O(\log{q})$. If each vertex of $G$ is adjacent to at least 
$\delta n$ vertices in each of the other vertex classes, $q$ is bounded and $n$ is large enough, then $G$ has a $K_q$-factor.
\end{abstract}

\section{Introduction}
In this paper we will consider simple graphs.
We mostly use standard notation: we denote by $V(F)$ and $E(F)$ the vertex
and the edge set of the graph $F$, $deg_F(x)$ is the degree of  
the vertex $x \in V(F)$ and $\delta(F)$ is the minimum
degree of $F$. 

Let $J$ be a fixed graph on $q$ vertices. If $q | |V(F)|$ and $F$ has a subgraph which consists of $|V(F)|/q$ 
vertex-disjoint copies of $J$, then we say that $F$ has a $J$-factor. 

\smallskip

A fundamental result in extremal graph theory is the following theorem of Hajnal and Szemer\'edi~\cite{HASZ70}:

\begin{theorem}[Hajnal and Szemer\'edi]
\label{hsz}
Let $G$ be a graph on $n$ vertices such that $\delta(G) \ge {q-1 \over q}n$. If $q$ divides $n$, then $G$ contains 
$n/q$ vertex-disjoint cliques of size $q$.
\end{theorem}

The theorem is obvious for $q=2$; the first non-trivial case $q=3$ was proved by K. Corr\'adi and A. Hajnal~\cite{COHA63}. 
The proof for arbitrary $q$ is notoriously hard, it was found by A. Hajnal and E. Szemer\'edi in 1970. 
  
\smallskip

We say that $F$ is {\it multipartite}, if its vertex set can be divided into classes which are independent 
sets. If the number of classes is $q$, then $F$ is $q$-partite. $F$ is a {\it balanced $q$-partite graph}, if these 
vertex classes are of the same size.
Let $F$ be a $q$-partite graph with vertex classes $A_1, A_2, \ldots, A_q$. We define the {\it proportional minimum degree} of $F$ by 

$$\widetilde{\delta}(F)=\min_{1 \le i \le q}{ \min_{v \in A_i} \{ {deg(v, A_j) \over |A_j|}: j \not =i \}}.$$

It is natural to investigate the multipartite version of Theorem~\ref{hsz}:

\begin{conjecture}\label{multi}
Let $G$ be a balanced $q$-partite graph on $qn$ vertices. There exists a positive constant $K$ such that if every vertex is adjacent to
at least ${q-1 \over q}n+K$ vertices in each of the other vertex classes, then  $G$ contains 
$n$ vertex-disjoint cliques of size $q$.
\end{conjecture}

Notice the extra additive constant: it turns out that it is necessary to have $K$ for odd $q$s. 
The conjecture is easily seen to hold for $q=2$. It was shown for $q=3$~\cite{MM} and $q=4$~\cite{MSz}. 
The proofs of  these latter cases are very involved.  
In this paper we show a relaxed version. For $k$ being a natural number let $h_k$ denote the $k$th harmonic number, that is, 
$h_k=1+{1 \over 2}+{1 \over 3}+ \ldots +{1 \over k}.$

\begin{theorem}\label{tetel}
Let $q\ge 3$ be an integer and $k_q=q-3/2+h_q/2.$ Then there exists an $n_0$ such that if 
$n>n_0$,  $G$ is a balanced $q$-partite graph on $qn$ vertices, and $\widetilde{\delta}(G)\ge 
{k_q \over k_q+1}$, then  $G$ has a $K_q$-factor.
\end{theorem} 

We also have the following corollary of Theorem~\ref{tetel}:

\begin{corollary}\label{kov}
Let $G$ be as above. Assume that $H$ is a fixed graph such that $\chi(H) \le q,$  and a constant number of vertex disjoint 
copies of $H$ is colored by $q$ colors such that we use every color and every color class has size $\kappa.$
If $\kappa$ divides $n,$
then $G$ has an $H$-factor. 
\end{corollary}

For proving Theorem~\ref{tetel} our main tools will be
the Regularity Lemma of Szemer\'edi~\cite{SZ76}, and the Blow-up Lemma~\cite{KSSz97, KSSz98}. We will
give a brief survey on 
the necessary notions in the second section.

\section{Main tools for the proof} \label{tools}

We introduce some more notation first. For any vertex $v$ of  the graph $G$, $deg_G(v,X)$ 
is  the number of  neighbors of $v$ in the set $X$, and $e(X,Y)$  is the number of  edges between  
the disjoint sets $X$ and $Y$.
$N_G(v)$ is the  set of neighbors of $v$ and $N_G(v,X)$  is the set of
neighbors of $v$ in $X$. For a set $S \subset V(G)$, $N(S)=\cup_{v \in S}N(v)$.

If $k$ is a natural number, and every vertex in the graph $G$ has degree $k$, then we call the graph $k$-regular,
or simply regular. On the other hand, for a real $\varepsilon \in (0,1)$  we will consider $\varepsilon$-regular pairs, 
these pairs play a crucial role in the Regularity Lemma of Szemer\'edi (more details follow later).

Let $F$ be a multipartite graph. Given certain vertex classes $A_{i_1}, \ldots, A_{i_s}$ we will denote the $s$-partite subgraph of $F$ spanned by these classes 
by $F(A_{i_1}, \ldots, A_{i_s}).$ Throughout the paper we will apply the relation ``$\ll$'': $a \ll b$, if $a$ is sufficiently smaller, than $b$.

\subsection{Factors of bipartite graphs}

Let $F$ be a bipartite graph with color classes $A$ and $B$. By the well-known K\"onig--Hall theorem there is a 
perfect matching in $F$ if and only if
$|N(S)| \ge |S|$ for every $S \subset A$. The following, while simple, is a very useful consequence of this result, we record
it here for future reference.

\begin{lemma} \label{egy-ketted}
If $F$ is a balanced bipartite graph on $2n$ vertices, and $deg(x) \ge n/2$ for every $x \in V(F)$, then there is a perfect 
matching in $F$.
\end{lemma}

If $f: V(F) \rightarrow {\rm {\bf  N}}$ is a function, then an $f$-factor is a subgraph $F'$ of $F$ such that $deg_{F'}(v)=f(v)$ for 
every $v \in V$. We will need special $f$-factors, namely when $f \equiv r$ for some $r \in N$. Then $F'$ is an 
$r$-regular subgraph  of $F$. If the minimum degree of $F$ is large enough, then one can find a sufficiently
dense spanning regular subgraph (see~\cite{Cs3}):

\begin{theorem}\label{regularis}
Let  $F(A,B)$ be a balanced bipartite graph on $2n$ vertices, and assume that
$\delta=\delta(F)/n \ge 1/2$. Then $G$ has a $\lfloor \rho(\delta) n \rfloor$--regular spanning subgraph, where 
$\rho(\delta) =  {\delta+ \sqrt{2 \delta-1} \over 2}.$
\end{theorem}

\subsection{Regularity Lemma}

The {\it density}  between disjoint sets $X$ and $Y$ is defined as:
$$d(X,Y)={{e(X,Y)}\over     {|X||Y|}}.$$
In   the   proof   of  Theorem~\ref{tetel},   Szemer\'edi's   Regularity
Lemma~\cite{SZ76, KS93} plays a pivotal  role. We  will need  the following
definition to state the Regularity Lemma.

\begin{definition}[Regularity condition] Let $\varepsilon >0$. A pair
  $(A,B)$      of     disjoint      vertex-sets     in      $G$     is
  \mbox{$\varepsilon$-regular}  if  for every  $X  \subset  A$ and  $Y
  \subset B$, satisfying
$$|X|>\varepsilon |A|,\ |Y|>\varepsilon |B|$$
we have
$$|d(X,Y)-d(A,B)|<\varepsilon.$$
\end{definition}
This  definition  implies  that   regular  pairs  are  highly  uniform
bipartite graphs; namely, the  density of any reasonably large subgraph
is almost the same as the density of the regular pair.

We will use the following form of the Regularity Lemma:

\begin{lemma}[Degree Form] \label{rl0} 
For every $\varepsilon>0$ there is an $M=M(\varepsilon)$ such
that if $G=(V,E)$ is any graph and $d\in [0,1]$ is any real number, then there is a partition of the
vertex set $V$ into $\ell+1$ clusters $W_0, W_1,\ldots,W_\ell$, and there is a subgraph $G'$ of $G$
with the following properties:  
\begin{itemize} 
\item $\ell \le M$, 
\item $|W_0|\le \varepsilon |V|$,
\item all clusters $W_i$, $i\ge 1$, are of the same size $m \ \left(\le \lfloor {|V|\over \ell}\rfloor<
\varepsilon |V| \right)$, 
\item $deg_{G'}(v)>deg_{G}(v)-(d+\varepsilon)|V|$ for all $v\in V$, 
\item $G'|_{W_i}=\emptyset$ ($W_i$ is an independent set in $G'$) for all $i\ge 1$, 
\item all pairs
$(W_i,W_j)$, $1\le i <j \le \ell$, are $\varepsilon$-regular, each with density either 0 or greater
than $d$ in $G'$.  
\end{itemize} 
\end{lemma}

\noindent Often  we  call  $W_0$  the  {\it exceptional  cluster}.
In the rest of the paper we will assume that $0< \varepsilon\ll  d \ll  1$.

\begin{definition}[Reduced graph]
Apply Lemma~\ref{rl0} to the graph $G=(V,E)$ with parameters $\varepsilon$ and $d$, 
and denote the clusters of the resulting partition by $W_0, W_1,\ldots, W_\ell$, $W_0$ being
the exceptional cluster. We construct a new graph $G_r$, the reduced graph of $G'$ in
the following way: The non-exceptional clusters of $G'$
are the vertices of the reduced graph $G_r$ (hence $|V(G_r)|=\ell$). We connect two vertices of
$G_r$ by an edge if the corresponding two clusters form an $\varepsilon$-regular
pair with density at least $d$.
\end{definition}

The following corollary is immediate:

\begin{corollary}\label{redukaltfok}
Apply Lemma~\ref{rl0} with parameters $\varepsilon$ and $d$ to the graph $G=(V,E)$ satisfying 
$\delta(G) \ge \gamma n$ $\ (|V|=n)$ for some $\gamma>0$. 
Denote $G_r$ the reduced graph of $G'$. Then $\delta(G_r) \ge (\gamma - \theta)\ell$,
where $\theta=2\varepsilon+d$.
\end{corollary} 

The lemma below states that the property of being balanced can be inherited by 
the reduced graph.

\begin{lemma} \label{balance}
Let $G$ be a balanced multipartite graph, then $G_r$ can be balanced as well.
\end{lemma}

\noindent {\bf Proof:} Trivial. \qedf

\medskip

Given an $\varepsilon$-regular pair $(A,B),$ we may increase $A$ and $B$ by adding some new vertices to
both. We expect that after this procedure the new pair will be $\eta$-regular for some small $\eta$, although $\eta > \varepsilon.$

\begin{lemma}\label{eps-reg-add}
Assume that $0 < \varepsilon \ll 1/K.$ Let $(A,B)$ be an $\varepsilon$-regular pair with $m=|A|=|B|,$ and add $K\varepsilon m$ 
vertices to $A$ and to $B.$ Then the resulting new pair is $2\sqrt{\varepsilon}$-regular.
\end{lemma}

\noindent {\bf Proof:} Simple computation. \hfill \qedf
\medskip

We will need the following simple lemma:

\begin{lemma} \label{felezo}
Let $(A,B)$ be an \mbox{$\varepsilon$-regular}--pair with density $d$ for some $\epsilon >0$.
We arbitrarily halve $A$ and $B$, getting the sets $A',A''$ and $B',B''$, respectively. Then the 
following holds: 
$(A',B')$ and $(A'',B'')$ are 
$2\varepsilon$--regular pairs with density at least $d-\varepsilon$.
\end{lemma}

\noindent {\bf Proof.} Trivial. \qedf

\medskip

A  stronger one-sided property of regular pairs is super-regularity:

\begin{definition}[Super-Regularity condition] Given a graph $G$ and
  two disjoint subsets  of its vertices $A$ and  $B$, the pair $(A,B)$
  is   \mbox{$(\varepsilon, \delta)${\rm    -super-regular}},   if   it   is
  $\varepsilon$-regular and furthermore,
$$deg(a)>\delta |B|,{\rm \ for\  all}\  a\in A,$$
and
$$deg(b)>\delta |A|,{\rm \ for\  all}\  b\in B.$$
\end{definition}

\medskip

Let $\varepsilon >0$ and assume that the pair $(A,B)$ is $\varepsilon$-regular with density $d.$ 
Mark those vertices of $A$ which have less than $(d-\varepsilon)|B|$ neighbors and those which have
more than $(d+\varepsilon)|B|.$ By the definition of $\varepsilon$-regularity, there can be at most 
$2\varepsilon|A|$ marked vertices in $A.$ Repeat the same procedure for $B$ so as to mark those
vertices which have too many or too few neighbors in $A.$ If we get rid of the marked vertices of $A$
and $B$ then we will have a \mbox{$(3\varepsilon, d-3\varepsilon)${\rm    -super-regular}} pair $(A', B')$.
That is, we proved that every regular pair contains a large super-regular pair:

\begin{lemma}\label{sub_super-reg}
Let $(A,B)$ be an $\varepsilon$-regular pair with density $d.$ Then it has a $(3\varepsilon, d-3\varepsilon)$-super-regular
subpair $(A', B')$ where $A'\subset A$, $|A'| = |A|-2\varepsilon|A|$ and  $B'\subset B$, $|B'| = |B|-2\varepsilon|B|.$
\end{lemma}

We will repeatedly make use of the following folklore result, which states that random subpairs
of $(\varepsilon, \delta)$-super-regular pairs are likely to be super-regular, with somewhat
weaker parameters:

\begin{proposition} \label{vagdalas}
Let $(A,B)$ be an $(\varepsilon, \delta)$-super-regular pair with density $d$ and $k$ be a positive integer. Assume that
$|A|=|B|=m$, and $k|m.$ Divide $A$ and $B$ into $k$ random subsets: 
$A=A_1 \cup A_2 \cup \ldots \cup A_{k}$ and $B=B_1 \cup B_2 \cup \ldots \cup B_k$, each having size
$m/k.$ Then with probability tending to one as $m$ tends to infinity we have that  $(A_i,B_j)$ is an 
$(\varepsilon', \delta')$-super-regular pair with  density $d'$ for every $1 \le i,j \le k$, where 
$\varepsilon' \le 2 \varepsilon$,  $\delta-\varepsilon \le \delta'$ and $d -\varepsilon \le d'.$ 
\end{proposition}

Let $G_r$ be the reduced graph of the graph $G$  such that edges in $G_r$ represent $\varepsilon$-regular 
pairs with density at least $d$.
Assume that $\widehat{G}_r$ is a cluster graph which we get by randomly splitting 
the clusters of $G_r$ into sub-clusters of equal size. The new sub-clusters will be called {\it split copies} of the original 
cluster, and we will use ``$\ \widehat{} \ $'' to indicate that we refer to a split copy.
 
Two split copies will be connected if they form an $\varepsilon'$-regular
pair with density $d'$ where $\varepsilon' \le 2 \varepsilon$ and $d' \ge d-\varepsilon.$ 
By the previous proposition if $W_iW_j \in E(G_r)$ and 
$\widehat{W}_i, \widehat{W}_j$ arose from $W_i$ and $W_j$ by the random splitting, then 
$\widehat{W}_i\widehat{W}_j \in E(\widehat{G}_r).$ We will call $\widehat{G}_r$
the {\it refinement} of $G_r.$

\subsection{Blow-up Lemma}

Let $H$ and $G$ be two graphs on $n$ vertices. Assume that we want to find an isomorphic copy of $H$ in $G$.
In order to achieve this one can apply a very powerful tool, the Blow-up Lemma of Koml\'os, S\'ark\"ozy and 
Szemer\'edi~\cite{KSSz97,KSSz98}.

\begin{theorem}[Blow-up Lemma]\label{blow-up}
Given a graph $R$ of order $r$ and positive integers $\delta, \Delta,$ there exists
a positive   $\varepsilon= \varepsilon(\delta ,\Delta ,r)$ such that the following holds: Let $n_1,n_2,\ldots,n_r$ be arbitrary positive 
parameters and let
us replace the vertices $v_1,v_2,\ldots,v_r$ of $R$ with pairwise disjoint sets $W_1,W_2, \ldots, W_r$ of sizes 
$n_1,n_2,\ldots,n_r$ (blowing up $R$). 
We construct two graphs on the same vertex set $V =\cup_i  W_i.$ The first graph $F$ is obtained by replacing each edge 
$v_iv_j \in E(R)$ with the complete bipartite graph between $W_i$ and $W_j.$ A sparser graph $G$ is constructed by replacing each edge 
$v_iv_j$ arbitrarily with an $(\varepsilon ,\delta )$-super-regular pair between $W_i$ and $W_j.$ If a graph $H$ with $\Delta (H)$   is embeddable 
into $F$ then it is already embeddable into $G.$
\end{theorem}

\section{The first stage of the embedding algorithm}

Since a $K_q$-factor is a subgraph, finding such a factor will be considered as an embedding problem. Let us denote the 
union of $n$ vertex-disjoint copies of $K_q$s by $H$. 
We will show Theorem~\ref{tetel} by exhibiting a randomized algorithm which with high probability will embed $H$ 
into $G$.  

The algorithm will proceed as follows: first, apply the Regularity Lemma to $G$ with appropriately chosen parameters
$0 < \varepsilon \ll d \ll 1,$ and get the balanced $q$-partite reduced graph $G_r$. The cluster classes of $G_r$ are denoted
by $A_1, A_2, \ldots, A_q$, here  $|A_1|=|A_2|=\ldots=|A_q|=\ell$. 

We will proceed in two stages. In the first stage we distribute the vertices of $H$ among the non-exceptional 
clusters of $G_r$: we require that (1) if two vertices of $H$ are adjacent, then they should be assigned to adjacent clusters of $G_r$; 
and (2) about the same number of vertices 
should be assigned to every cluster -- the difference cannot be larger than $o(n)$. These requirements will
be achieved via finding a $K_q$-factor in a cluster graph $\widehat{G}_r,$ which is a refinement of $G_r.$
The recursive algorithm to construct the clique-factor in $\widehat{G}_r$ is called the {\it Factor Finder Algorithm.}

Having the above mentioned clique-factor we take any surjective function $\phi$ which assigns $q$-cliques of $H$ to $q$-cliques of the
factor in $\widehat{G}_r$ in such a way that $|\phi^{-1}(C)|=|\phi^{-1}(C')|$ for every two cliques $C, C'.$  Such a function is obvious to find.
Then we assign 
the vertices of the clique $\widetilde{C} \in H$ to the clusters of $\phi(\widetilde{C})$ in the obvious way.
This assignment is easily seen to satisfy (1) and (2). 

We will finish the embedding in the second stage by the help of the Blow-up Lemma. This is a technically somewhat challenging
part, however, this stage is more routine.

The Factor Finder algorithm is a recursive algorithm, with base case $q=2.$ For an easier understanding we will consider 
the case $q=3$ in greater detail, and then generalize the method for larger $q$s.

\subsection{The Factor Finder Algorithm}

Given the graph $G$, we apply the Degree Form of the Regularity Lemma with parameters $\varepsilon$ and $d$ such 
that $0 < \varepsilon \ll d \ll 1$.
Then we find the reduced graph $G_r$. We assume that it is a balanced $q$-partite reduced graph on $q \ell$ vertices with 
$\widetilde{\delta}(G_r) \ge k_q /(k_q+1)$ where $k_q=q-3/2+h_{q-1}/2$. Observe that this is not necessarily the case, since
$\widetilde{\delta}(G) \ge k_q /(k_q+1)$ does not imply the above bound for $\widetilde{\delta}(G_r),$ in general we may lose some edges
when discarding the irregular pairs and also when putting vertices to $W_0$ (recall Corollary~\ref{redukaltfok}). It turns out that with 
a proportional 
minimum degree this large we will have room to spare in case $q \ge 3.$ Therefore, the algorithm will find the $K_q$-factor even in case
$\widetilde{\delta}(G_r) =k_q/(k_q+1)-\gamma_q,$ where $\gamma_q$ is a function of $q.$ We will discuss the details at the end.

In what follows we will denote $k_q/(k_q+1)$ by $\widetilde{\delta},$ 
and the cluster classes of $G_r$ will be denoted by $A_1, A_2, \ldots, A_q.$
Recall, that our goal
is to show that $H \subset G,$ where $H$ is the disjoint union of $n$ copies of $K_q$s.

\medskip

\noindent {\bf The first case: $q=2$}

First, notice that in this case $k_q=1,$ therefore, $\widetilde{\delta}=1/2.$ It is straightforward to find a $K_2$-factor 
(a perfect matching) in a balanced bipartite graph $G_r$ with a proportional minimum degree this large 
(Lemma~\ref{egy-ketted}). 


\smallskip

Assume, that the proportional minimum degree is a bit larger, it is $1/2+\psi$ for
some $0<\psi<1.$ Then one can introduce some randomness in finding the
perfect matching. Pick $\psi\ell/2$ clusters randomly from the first vertex class, and
find neighbors for them randomly. Then pick $\psi\ell/2$ clusters randomly from the 
other vertex class, and find neighbors for them randomly. This way we have found 
random neighbors for $\psi \ell$ clusters. In the leftover the minimum degree is
sufficiently large for having a perfect matching. Therefore, we can find a perfect
matching in such a way that $\psi \ell$ clusters have randomly chosen neighbors.

As it turns out later on, this small extra randomness will be very helpful. When finishing
the embedding of $H$ we need
a bit larger proportional minimum degree, than $1/2$ at the end in order to perform this
procedure (recall, that we use
recursion), but that will be provided for $q \ge 3.$

\medskip

\noindent {\bf Finding a triangle factor}
\smallskip

As a warm-up we discuss this case in details. 
First, apply Theorem~\ref{regularis} for the graphs $G_r(A_1,A_2)$ and $G_r(A_1,A_3).$ We get two $\mu$-regular bipartite 
graphs $R(A_1,A_2)$ and $R(A_1,A_3),$ with $\mu = \rho(\widetilde{\delta}) \ell.$ (Since we can delete as many 1-factors from a 
regular bipartite graph 
as we please and still get a regular bipartite graph, we may assume that  $\mu=\rho(\widetilde{\delta})\ell$.) We let 
$R$ to be a $3$-partite graph on $A_1 \cup A_2 \cup A_3$ such that $E(R)=E(R(A_1, A_2)) \cup E(R(A_1, A_3)).$ 
It is easy to see, that  $deg_R(b)=\mu$ for every $b \in A_2 \cup A_3.$ 

We are going to cut the clusters
of $A_2 \cup A_3$ randomly into $\mu$ sub-clusters of equal size. The new cluster classes are denoted by $\widehat{A}_2$ and 
$\widehat{A}_3.$ 
Roughly speaking, we will assign the split copies of $\widehat{A}_2 \cup \widehat{A}_3$ to the clusters of $A_1,$ such that every cluster of $A_1$ will receive $2\mu$ split
copies, and every split copy will be assigned to exactly one cluster in $A_1.$

More formally, let us define a surjective function $\sigma$: its domain is the set of split copies, and its range is $A_1.$ 
It satisfies the following requirements: whenever 
$U \in A_2 \cup A_3,$ and $\widehat{U}$ is a split copy of $U$, then $\sigma(\widehat{U}) \in N_R(U),$ 
moreover, if $\widehat{U}$ and $\widehat{U}'$ are
different split copies of $U$, then $\sigma(U) \neq \sigma(\widehat{U}').$ 
For every $W \in A_1$ we introduce two sets associated with it: $$N_2(W)=\{\widehat{U}: \widehat{U} \in \widehat{A}_2, \ 
\sigma(\widehat{U})=W\},$$ and   
$$N_3(W)=\{\widehat{U}: \widehat{U} \in \widehat{A}_3, \ \sigma(\widehat{U})=W\}.$$ It is easy to see, that every cluster of 
$\widehat{A}_2 \cup \widehat{A}_3$ will participate in
one of the $N_i(W)$ sets, and $|N_i(W)|=\mu$ for $i=2, 3$ and every $W \in A_1.$

Our next goal is to show, that $\widehat{G}_r(N_2(W), N_3(W)),$ the induced subgraph of the refinement of $G_r$ on
$N_2(W)$ and $N_3(W)$ has a perfect matching $M(W)$ for every $W \in A_1.$ Having this perfect matching at hand
we can construct $\mu$ triangles for every $W \in A_1$: cut the clusters of $A_1$ randomly into $\mu$ sub-clusters,
and assign the split copies of $W$ to the edges of $M(W)$ bijectively. This way we construct triangles each having cluster
size $m/\mu.$

Hence, what is left: for every $W \in A_1$ find the perfect matchings in the bipartite subgraphs $\widehat{G}_r(N_2(W), N_3(W)).$ 
We claim that the minimum degree in these bipartite graphs is in fact sufficiently 
large to guarantee the existence of a perfect matching in it. For that we will show that every cluster is adjacent to at least half of the 
clusters in the other class. 
We use a simple claim which we record here for future purposes.

\begin{claim} \label{fokarany}
Let $F=(V,E)$ be a graph and let $S \subset V$. Then every $u \in V$ is adjacent to at least 
a ${\delta(F)-(|V|-|S|) \over |S|}$ proportion of the vertices of $S$.
\end{claim}

\noindent {\bf Proof:} Obvious. \qedf

\smallskip

Now let $\widehat{U} \in N_2(W)$ be an arbitrary cluster. By Claim~\ref{fokarany} $\widehat{U}$ is adjacent to at least 
$(\widetilde{\delta} - (1-\mu/\ell))\ell/\mu$ 
proportion of the vertices of $N_3(W).$ Similarly, every $\widehat{U} \in N_3(W)$ is adjacent to at least 
$(\widetilde{\delta} - (1-\mu/\ell))\ell/\mu$
proportion of $N_2(W).$
Easy calculation shows that if $\widetilde{\delta} = 0.68$, then 
$\mu/\ell=\rho(\widetilde{\delta})=0.64$, and  
$$\widetilde{\delta}(\widehat{G}_r(N_2(W), N_3(W))) \ge {(\widetilde{\delta} - (1-\mu/\ell))\ell \over \mu}=0.5.$$ 
This implies the existence of a perfect matching in
$\widehat{G}_r(N_2(W), N_3(W)),$ hence, as we discussed above, this can be extended into a triangle factor in
the refinement $\widehat{G}_r.$ Notice, that $0.68 < 0.6923<{3-3/2+h_2/2 \over 3-3/2+h_2/2 +1}$, i.e., we have
found a triangle factor in $\widehat{G}_r$ with a smaller bound that is required by Theorem~\ref{tetel}. 

This latter fact will be important for us later on. Recall the discussion
of case $q=2$. Obviously, for $k_3=3/2+h_2/2$ the proportional minimum degree 
will be larger than $1/2$ when it comes to finding the perfect matchings in the $\widehat{G}_r(N_2(W), N_3(W))$ graphs. Hence,
we can perform the randomized procedure for finding the perfect matchings.

As we noted above, having a triangle factor in $\widehat{G}_r$ allows us to find the good pre-assignment easily.
We remark, that the cluster size in $\widehat{G}_r$ is $m \over \mu$, and the number of clusters is 
$\mu \ell=\rho(\widetilde{\delta})\ell^2.$

\medskip

\noindent {\bf The general case}

Assume now that $q>3$. We will apply induction on $q,$ and assume that if the proportional minimum degree
in a balanced $(q-1)$-partite cluster graph $F$ is at least $k_{q-1}/(k_{q-1}+1),$ then $F$ has a $K_{q-1}$-factor.

We are given $G_r$, a balanced $q$-partite graph with vertex classes $A_1, A_2, \ldots, A_q$ such that 
$\widetilde{\delta}(G_r) \ge k_q/(k_q+1).$
This time our goal will be to find a $K_q$-factor in a refinement $\widehat{G}_r.$ 

Set $\mu=\rho(\widetilde{\delta}(G_r))\ell$. We consider the bipartite subgraphs $G_r(A_1, A_i)$  
and apply Theorem~\ref{regularis} to get 
the $\mu$-regular bipartite graphs $R(A_1, A_i)$ for every $2 \le i \le q.$ Let $R$ be a $q$-partite graph such 
that $V(R)=\cup_{i \ge 1}A_i$ and $E(R)=\cup_{i \ge 2}E(R(A_1, A_i)).$ As before, $deg_R(U)=\mu$ where
$U \in A_2 \cup \ldots \cup A_q.$

Similarly to the case $q=3$ we randomly split every cluster in $A_2 \cup A_3 \cup \ldots \cup A_q$
into $\mu$ sub-clusters of equal size thereby getting $\widehat{A}_i$ from $A_i$ for $2 \le i \le q.$

We define a surjective function $\sigma$: its domain is the set of split copies, and its range is $A_1.$ It satisfies the
following requirements: whenever 
$U \in A_2 \cup \ldots \cup A_q,$ and $\widehat{U}$ is a split copy of $U$, then $\sigma(\widehat{U}) \in N_R(U),$ moreover, 
if $\widehat{U}$ and $\widehat{U}'$ are
different split copies of $U$, then $\sigma(\widehat{U}) \neq \sigma(\widehat{U}').$ 
For every $W \in A_1$ we introduce $q-1$ sets associated with it: $$N_i(W)=\{\widehat{U}: \widehat{U} \in \widehat{A}_i, \ 
\sigma(\widehat{U})=W\}$$ 
for $2 \le i \le q.$
It is easy to see, that every cluster of $\widehat{A}_2 \cup \ldots \cup \widehat{A}_q$ will participate in
one of the $N_i(W)$ sets, and $|N_i(W)|=\mu$ for every $2 \le i \le q$ and every $W \in A_1.$

Let us consider the balanced $(q-1)$-partite graphs $\widehat{G}_r(N_2(W), \ldots, N_q(W))$ for every $W \in A_1.$ 
As before, we can lower bound the proportional 
minimum degree in these graphs by the help of Claim~\ref{fokarany}:
 $$\widetilde{\delta}(\widehat{G}_r(N_2(W), \ldots, N_q(W)))\ge {(\widetilde{\delta} - (1-\mu_q/\ell))\ell \over \mu_q}.$$ 
In case $q=3$ we had to check whether this quantity was at least $1/2,$ this time we have to check that this number is sufficiently 
large so as to guarantee the existence of a $K_{q-1}$-factor in these graphs. 

Say, that we can find a $K_{q-1}$-factor $M(W)$ for every $W \in A_1.$
Then we construct the desired $K_q$-factor in the following way:
cut the clusters of $A_1$ randomly into $\mu$ sub-clusters,
and assign the split copies of $W$ to the $(q-1)$-cliques of $M(W)$ bijectively. This way we get $\mu$ cliques of size $q$ each 
having clusters of size $m/\mu.$

In Lemma~\ref{induction} below we will prove that $\widetilde{\delta} \ge {k_q \over k_q+1}$ is sufficiently large, that is, 
the proportional minimum degree
in  $\widetilde{\delta}(\widehat{G}_r(N_2(W), \ldots, N_q(W)))$ is at least ${k_{q-1} \over k_{q-1}+1},$ which, by the 
induction hypothesis implies the existence of a $K_{q-1}$-factor in $G_r^j(A_2^j, A_3^j, \ldots, A_q^j)$ for every $1 \le j \le \ell.$ 

\begin{lemma} \label{induction}
$$\widetilde{\delta}(\widehat{G}_r(N_2(W), \ldots, N_q(W))) \ge {k_{q-1} \over k_{q-1}+1}$$ for every $1 \le j \le \ell$ if $q\ge 3.$
\end{lemma}

\noindent {\bf Proof:} By Claim~\ref{fokarany} we have that 
$$\widetilde{\delta}(\widehat{G}_r(N_2(W), \ldots, N_q(W))) \ge {{k_{q} \over k_{q}+1}-(1-\rho({k_{q} \over k_{q}+1})) 
\over \rho({k_{q} \over k_{q}+1})}.$$
Since $$\rho({k_{q} \over k_{q}+1})={{k_q \over k_q+1}+\sqrt{{k_q-1 \over k_q+1}} \over 2},$$
we get the following lower bound for the proportional minimum degree in $\widehat{G}_r(N_2(a), \ldots, N_q(a))$:
$$2{{k_q-2 \over 2(k_q+1)}+{1 \over 2}\sqrt{k_q-1 \over k_q+1} \over {k_q \over k_q+1}+\sqrt{k_q-1 \over k_q+1}}=
1-{{2 \over k_q+1} \over {k_q \over k_q+1}+\sqrt{k_q-1 \over k_q+1}}=1 -{2 \over k_q+\sqrt{(k_q-1)(k_q+1)}}.$$

We will show that $$1 -{2 \over k_q+\sqrt{(k_q-1)(k_q+1)}} > {k_{q-1} \over k_{q-1}+1}=1-{1 \over k_{q-1}+1}$$ 
is a valid inequality. Equivalently, we claim that 
$$1 -{2 \over k_q+\sqrt{(k_q-1)(k_q+1)}} -(1-{1 \over k_{q-1}+1})={1 \over k_{q-1}+1}-{2 \over k_q+\sqrt{(k_q-1)(k_q+1)}}>0.$$
This is implied by
$$k_q+k_q\sqrt{1 -{1 \over k_q^2}} > 2(k_{q-1}+1).$$ Since $k_q=k_{q-1}+1+1/(2q-2),$ the above is a consequence of
$$k_q\sqrt{1-{1 \over k_q^2}} > k_q+{1 \over q-1},$$ which is easily seen to hold for $q \ge 3.$ \hfill \qedf

\medskip

We have proved that the Factor Finder algorithm can construct a $K_q$-factor in $\widehat{G}_r$ provided that 
$\widetilde{\delta}(G_r) \ge k_q/(k_q+1).$ Observe that apart from the case $q=2$ we have room to spare in the
proportional minimum degree. That is, the algorithm will complete its task successfully even in case 
$\widetilde{\delta}(G_r) = k_q/(k_q+1) -\gamma_q$ if $q \ge 3$ and $\gamma_q$ is sufficiently small. Let us choose $\varepsilon$
and $d$ such that $0 < \varepsilon \ll d \ll \gamma_q,$ and apply the Regularity Lemma. With this choice, by Corollary~\ref{redukaltfok},
$\widetilde{\delta}(G_r) \ge k_q/(k_q+1) -\gamma_q.$ Hence, if $\varepsilon$ and $d$ are sufficiently small then the Factor Finder
algorithm will find the clique factor.

We remark that the bound of $k_q=q-3/2+h_{q-1}/2$ could be improved somewhat. We didn't want to 
optimize on this bound. It already gives the correct order of magnitude for our embedding method:
$k_q=q+O(\log{q}),$ without having tedious computations in the proof of the lemma.

\medskip

\noindent {\bf More on the Factor Finder algorithm}

\smallskip

Let us explore more properties of the Factor Finder algorithm, which will be useful later on.
Set $s_1(q)=\ell$ for every $q \ge 3.$ Given a cluster $W \in A_1$ we denote its degree in $R(A_1, A_i)$ by $s_2(q)$, that is, 
$s_2(q)=\rho(\widetilde{\delta}(G_r))\ell.$ 
The recursive process guarantees that we can construct a $K_{q-1}$-factor in the $q-1$ neighborhoods of $W$, each having
size $s_2(q).$ Now for finding the $K_{q-1}$-factor we again apply recursion, and want to find a $K_{q-2}$ factor in $s_2(q)$ different
balanced $(q-2)$-partite graphs. The size of the vertex classes of these balanced graphs will be denoted by $s_3(q).$ In general, when
proceeding with the recursion, step-by-step we construct balanced $(q-i)$-partite 
graphs, in which we look for
a $K_{q-i}$-factor. The number of these graphs is $s_1(q)\cdot s_2(q) \cdots s_{i}(q).$ The number of clusters in a class of these balanced graphs are denoted by $s_{i}(q).$
We stop $i=q-1,$ when we arrive to balanced bipartite graphs, in which we are looking for perfect matchings.

We can compute the number of cliques in the $K_q$-factor which contain some split copy of a given cluster.

\begin{lemma} \label{klikk-fok}
Let $U$ be an arbitrary cluster in $G_r.$ The split copies of $U$ appear in
$\Pi_{i=2}^{q-1}s_{i}(q)$ cliques in the $K_q$-factor of $\widehat{G}_r.$ 
\end{lemma}

\noindent {\bf Proof:}
We want to apply induction, but for doing that we have to be careful. The statement we will prove by induction is as
follows: 

\noindent {\bf Claim:} Let $F$ be a balanced $a$-partite cluster graph with cluster classes of size $\ell,$ and $W$ be a cluster of $F.$ 
If $\widetilde{\delta}(F)\ge k_j/(k_j+1)$ where $j \ge a,$ and we apply the Factor Finder algorithm then the number of $a$-cliques 
containing a split copy of $W$ is 
$\Pi_{i=2}^{a-1}s_{i}(j).$ 

It is easy to see
that this statement is stronger than that of the lemma. Notice, that we have to keep track of the size of the cluster classes, too.

We show that in case $a=3$ the above statement holds. Let $j \ge 3.$ 
First assume that $U \in A_1.$ The algorithm finds the neighborhoods $N_2(U) \subset A_2$ and $N_3(U) \subset A_3,$
both having size $s_2(j).$ Next we look for a perfect matching between these two sets, every edge of this matching
with $U$ will result in a triangle. Hence, the number of triangles having a split copy of $U$ is $s_2(j).$ 

Suppose, that $U \in A_2,$ and let $ W \in N_R(U, A_1)$ be arbitrary. Then there will be triangle
which contains a split copy of $W$ and a split copy of $U$. Since this holds for every cluster of $N_R(U, A_1),$ and this set
has $s_2(j)$ clusters, there are $s_2(j)$ triangles which contain a split copy of $U.$

Assume now that $a >3$ and that the induction hypothesis holds up to $a-1.$ Let $j \ge a.$
As above, we begin with the case $U \in A_1.$ The algorithm
first finds an $(a-1)$-partite cluster graph in which every cluster class has size $s_2(j),$ and $U$ is adjacent to every cluster of this
graph. We want to find a $K_{a-1}$-factor in some refinement of it by the Factor Finder algorithm. Let $W$ be an arbitrary cluster 
from the ``first" cluster
class of the $a-1$ classes. We have $s_2(j)$ possible choices for $W.$ 
The following is easy to see: for $1\le i \le a-2$ the cluster classes of the $(a-i)$-partite graphs constructed by the Factor Finder 
algorithm will be of size $s_{i+1}(j).$ 
Hence, applying the induction hypothesis, there are $\Pi_{i=2}^{a-2}s_{i+1}(j)$ cliques on $a-1$ clusters which contain a split copy of 
$W.$
We have $s_2(j)$ choices for $W,$ therefore, the number of $q$-cliques containing a split copy of $U$ is 
$s_2(j)\Pi_{i=2}^{a-2}s_{i+1}(j)=\Pi_{i=2}^{a-1}s_{i}(j).$ 

\smallskip

Finally, we consider the case $a>3$ when $U \in A_t$ for $t>1.$
In the first step there are $s_2(j)$ clusters of $A_1$ such that these are adjacent to $U$ in $R(A_1, A_t).$ Let $W$ be any of these clusters.
Consider the $(a-1)$-partite cluster graph which is constructed for $W$ by the algorithm. This cluster graph has classes of size $s_2(j).$
As above, we can apply induction, and get that the algorithm finds $\Pi_{i=2}^{a-2}s_{i+1}(j)$ cliques on $a-1$ clusters which contain a 
split copy of $U.$ We repeat this for every cluster in $N_R(U, A_1),$ that results in $s_2(j)$ different $(a-1)$-partite graphs. In each of these we find $\Pi_{i=2}^{a-2}s_{i+1}(j)$ cliques on $a-1$ clusters containing a split copy of $U.$ Overall, split copies of $U$ appear in 
$s_2(j)\Pi_{i=2}^{a-2}s_{i+1}(j)=\Pi_{i=2}^{a-1}s_{i}(j)$ cliques on $a$ clusters. 
\hfill \qedf

\smallskip

Obviously, $s_1(q) > s_2(q) > s_3(q) > \ldots s_{q-1}(q) > {2\ell \over k_q+1}$ 
for $q \ge 3.$
The last inequality follows from Claim~\ref{fokarany} and the fact that the proportional minimum degree in the last graph is $\ge 1/2+ \psi_q$ for some positive
constant $\psi_q$ depending only on $q.$ (Recall that $k_3/(k_3+1)-0.68 >0.01$, hence, $\psi_3 >0.01$, and because
of Lemma~\ref{induction} the property of $\psi_q$ being positive is inherited for larger values of $q.$) 
Observe, that the overall number of cliques in the $K_q$-factor is $\Pi_{i=1}^{q-1}s_{i}(q)=\nu_q \ell^{q-1}$, where $\nu_q$ is a constant.
This implies, that
the cluster size in $G_r$ is $m=\Pi_{i=2}^{q-1}s_{i}(q) \widehat{m}$, where $\widehat{m}$ is the 
common cluster size in the refinement $\widehat{G}_r,$ and the number of clusters in $\widehat{G}_r$ is 
$\widehat{\ell}=\Pi_{i=1}^{q-1}s_{i}(q)=\nu_q \ell^{q-1}.$

\section{Second stage -- Finishing the proof of Theorem~\ref{tetel}}

In this section we discuss how to finish the embedding of $H$ into $G.$ Observe, that by applying
Lemma~\ref{sub_super-reg}, Proposition~\ref{vagdalas} and the Blow-up Lemma we are able to embed most of $H$
into $G$: The edges in the cliques of the $K_q$ factor of $\widehat{G}_r$ represent $\varepsilon'$-regular pairs, which by 
Lemma~\ref{sub_super-reg} can be made super regular. Applying the Blow-up lemma we get that most of $H$ can be embedded 
into $G$, at most $3\varepsilon' n+|W_0|\le 4\varepsilon' n$ vertices are left out, here $\varepsilon'$ is a constant multiple of
$\varepsilon.$ 

Our main goal in this section is to embed the {\it whole} of $H$ by the help of the Blow-up Lemma. For that we will try to find a 
$K_q$-factor in such a way that 
every edge in the cliques will represent $(\eta, d-\eta)$-super-regular pairs, where $\varepsilon$ will be a function of $\eta.$ Moreover, 
every vertex of $G$ will sit in a cluster of some clique, and every cluster will have the same size.
We will achieve this goal in a few steps. First we discard those vertices from the cliques which do not have many neighbors in other
clusters of the cliques, and put them to $W_0,$ the exceptional cluster. 
Secondly, we will distribute the vertices of $W_0$ such that every edge in the $K_q$-factor will 
represent a super-regular
pair. Finally, we move vertices between clusters so as to get equal size clusters in the cliques, but keep super-regularity, we
call this the {\it balancing step}. Then we will
apply the Blow-up Lemma.  
 
We need an important lemma, which will be crucial for making the cluster
sizes equal in every clique. In order to state it, let us 
define $q$ directed graphs: $L_1, L_2, \ldots, L_q.$ Here $V(L_i)=\widehat{A}_i,$ the class containing the split copies of
the clusters of the $i$th class.
Let $\widehat{U}_1, \widehat{U}_2  \in \widehat{A}_i,$ we will have the directed edge 
$(\widehat{U}_1,\widehat{U}_2) \in E(L_i)$, if $\widehat{U}_1$ is adjacent to all the clusters of the $q$-clique
which contains $\widehat{U}_2$ except $\widehat{U}_2$ itself. That is, if $\widehat{W}$ is a cluster of this clique, then the 
$(\widehat{U}_1,\widehat{W})$ pair is $\varepsilon'$-regular. We will
also say that $\widehat{U}_1$ is adjacent to the clique of $\widehat{U}_2.$
We will show the following: 

\begin{lemma}\label{szomszedsag}
Let $U_1, U_2 \in A_i$ for some $1 \le i \le q,$ and let $\widehat{U}_1$ be any split copy of $U_1$ in $\widehat{G}_r.$ Then with probability at least $1 - 1/(2q^2\ell^2)$ 
there are more than ${1 \over 8} \psi^2_qs_{q-1}(q)\Pi_{i=3}^{q-1}s_i(q)$ split copies of $\widehat{U}_2$ such that 
$\widehat{U}_1$ is adjacent to its clique.
\end{lemma}

The main message of Lemma~\ref{szomszedsag} is that 
out of the $\Pi_{i=2}^{q-1}s_i(q)$
cliques in the factor which contain some split copy of $U_2$ a constant proportion is adjacent to some split copy of $U_1,$
independently of the choice of $U_1$ and $U_2.$

\smallskip

\noindent{\bf Proof:} 
We will follow the line of arguments of the proof of Lemma~\ref{klikk-fok}. The extra cluster $U_1$ can be considered
as having one more cluster class. 
More precisely, the effect of having $U_1$ is as follows. When computing the number of cliques having
a split copy $\widehat{U}_2,$ at every step we have to take into account whether the clusters are in the neighborhood 
of $U_1.$ This shrinks the sizes: if the cluster class size in question is $s_i(q),$ then out of this many clusters at least 
$s_{i+1}(q)$ is adjacent to $U_1.$ 

This estimation works smoothly until at the end we have to find a
perfect matching in a bipartite graph having cluster classes of size $s_{q-1}(q)$ each. Then $U_1$ is adjacent to at least
$(1/2+\psi_q)s_{q-1}(q)$ clusters in both classes. Recall that we find the perfect matching in the following way: 
We randomly, independently, with probability $\psi_q/2$ choose clusters,
and pick a random vacant neighbor for those. The rest can get a neighbor by any algorithm for finding a perfect matching.
Suppose that we choose $W$ for having a random neighbor. Since $W$ and $U_1$ have a common neighborhood of size
at least $\psi_q s_{q-1}(q),$ the probability that $W$ will get a neighbor in the
perfect matching which is adjacent to $U_1$ is at least $\psi_q.$ The expected number of cliques containing some split copy 
$\widehat{U}_2$ and being
adjacent to $U_1$ is at least ${1 \over 4} \psi^2_qs_{q-1}(q)\Pi_{i=3}^{q-1}s_i(q),$ here we applied the bound of 
Lemma~\ref{klikk-fok}. Standard probabilistic reasoning -- use e.g.,
Azuma's inequality -- shows that $U_1$ will be adjacent to at least ${1 \over 8} \psi^2_qs_{q-1}(q)\Pi_{i=3}^{q-1}s_i(q)$ 
split copies of $U_2$ with probability at least $1-1/(2q^2\ell^2).$ \hfill\qedf

\medskip

Observe, that if $\widehat{U}, \widehat{U}'$ are split copies of $U,$ then $\widehat{U}$ is adjacent to the clique of $\widehat{U}'.$
Together with the so called union bound in probability theory this implies the following:

\begin{corollary}\label{L-path}
With positive probability there are at least ${1 \over 8} \psi^2_qs_{q-1}(q)\Pi_{i=3}^{q-1}s_i(q)$ vertex disjoint directed paths of length
at most two between any two clusters in $L_i,$ for every $1 \le i \le q.$
\end{corollary}

\medskip

We have acquired the knowledge to achieve our main goal, in the rest of the section we discuss how to finish the embedding 
step by step. 

In the first step we make every edge in the cliques of the factor super-regular by applying Lemma~\ref{sub_super-reg}, the discarded
vertices will be put to $W_0.$ 
Then the enlarged extremal cluster $W_0$ will be larger, but still remain reasonably small: $|W_0| \le \varepsilon' n,$
where $\varepsilon'$ is a constant multiple of $\varepsilon.$

In the second step we will distribute the vertices of $W_0$ among the $\widehat{\ell}$ clusters of $\widehat{G}_r.$ Let $v \in W_0$ and 
$\widehat{U}$ be a cluster. We say that $v$ is
adjacent to the clique of $\widehat{U}$ if $v$ has at least $d\widehat{m}$ neighbors in every cluster in the clique of $\widehat{U},$ except
in $\widehat{U}$ itself. Notice, that the proof of Lemma~\ref{szomszedsag} shows, that for every $v \in W_0$ there are at least 
${1 \over 8} \psi^2_q\ell s_{q-1}(q)\Pi_{i=3}^{q-1}s_i(q)$ clusters such that $v$ is adjacent to their cliques. Since the number of cliques is
$\widehat{\ell}=\Pi_{i=1}^{q-1}s_{i}(q),$ every vertex is adjacent to $c_q \widehat{\ell}$ cliques, where 
$c_q=s_{q-1}(q)\psi_q^2/(8s_2(q)).$ 
 
When distributing the vertices of $W_0$ we are allowed to put a vertex $v$ to a cluster $\widehat{U}$ if $v$ is adjacent to the clique 
of $\widehat{U}.$ We pay attention to distribute the vertices {\it evenly}, that is, at the end no cluster will get more than 
$|W_0|/(c_q\widehat{\ell})$ new vertices from $W_0.$ Since every vertex is adjacent to many cliques, this can be achieved.
After this step every edge of every clique in the $K_q$-factor will represent super-regular pairs. 

It is possible, that
the clusters have different sizes in a clique, hence, we have to perform the balancing algorithm. 
For that we partition the clusters of $\widehat{G}_r$ into three sets: $S_<, S_=$ and $S_>.$ $S_<$ contains those clusters which have
less than $n/\widehat{\ell}$ vertices, $S_>$ contains those clusters which have more than $n/\widehat{\ell}$ vertices, and $S_=$
contains the rest with equality. We will apply Corollary~\ref{L-path} in order to find directed paths from clusters of $S_>$ to clusters in $S_<.$

Say, that $\widehat{U}_1 \in S_>, \widehat{U}_2 \in S_<$ and there is a path of length one between them, that is, 
$\widehat{U}_1\widehat{U}_2 \in E(L_i)$ for some $1 \le i \le q.$
Then the vast majority of the vertices of $\widehat{U}_1$ are adjacent to the clique of $\widehat{U}_2.$ Pick as many as
needed (and possible) among these and place them to $\widehat{U}_2.$ 
If the path is of length two, then choose a cluster $\widehat{U}_3$
such that $\widehat{U}_1\widehat{U}_3$ and $\widehat{U}_3\widehat{U}_2$ belong to $E(L_i).$ Again, the vast majority of the
vertices in $\widehat{U}_1$ are adjacent to the clique of $\widehat{U}_3$ and the vast majority of the vertices
of $\widehat{U}_3$ are adjacent to the clique of $\widehat{U}_2.$ Hence, by placing vertices from $\widehat{U}_1$ to 
$\widehat{U}_3$ and the same number of vertices from $\widehat{U}_3$ to $\widehat{U}_2$ we can decrease the
discrepancy of $\widehat{U}_1$ and $\widehat{U}_2$ such that we keep the edges super-regular in all the cliques in question.
Observe, that we can perform the balancing algorithm such that we do not take out more than $|W_0|/(c_q \widehat{\ell})$ 
vertices from any of the clusters, and do not put in more than $|W_0|/(c_q \widehat{\ell})$ vertices to any of the cluster. 

We can apply Lemma~\ref{eps-reg-add}, and get that the edges of the cliques represent 
$(\widehat{\varepsilon}, \widehat{d})$-super-regular
pairs, where $\widehat{\varepsilon} \le C\sqrt{\varepsilon}$ and $\widehat{d} \ge d-\widehat{\varepsilon},$ and $C$ is a constant.
 
At this point we can recognize, that with positive probability all conditions of Lemma~\ref{blow-up}
are satisfied if $\varepsilon$ is sufficiently small and $\varepsilon \ll d \ll 1.$ From this the proof of Theorem~\ref{tetel} follows.

\bigskip

\noindent {\bf Proof of Corollary~\ref{kov}:}
We can embed vertex disjoint copies of $H$ as follows:
first, find a $K_q$-factor in $G.$ Then color some vertex disjoint union of copies of $H$ by $q$ colors such
that every color class has size $\kappa$ and every color is used. Call this colored graph $\widetilde{H}.$
It is easy to see that $G$ has an $\widetilde{H}$-factor: we embed the copies of $\widetilde{H}$ in the 
cliques of the $K_q$-factor by the help of the Blow-up Lemma. \hfill \qedf

\bigskip

\noindent {\bf Acknowledgment}
The author would like to thank P\'eter Hajnal and Endre Szemer\'edi  for the helpful conversations.


\begin{thebibliography}{20}

\bibitem{COHA63} H. Corr\'adi and A. Hajnal (1963) On the Maximal Number of
Independent Circuits in a Graph, Acta Math. Hung., {\bf 14}, 423-439.

\bibitem{Cs3} B. Csaba (2007) Regular Spanning Subgraphs of Bipartite Graphs of High Minimum Degree, The Electronic Journal of Combinatorics, \#N21.  

\bibitem{HASZ70} A. Hajnal and E. Szemer\'edi (1970) Proof of a Conjecture
of Erd\H{o}s, in ``Combinatorial Theory and Its Applications, II''
(P. Erd\H{o}s, and V. T.  S\'os, Eds.), Colloquia Mathematica Societatis
J\'anos Bolyai, North-Holland, Amsterdam/London.

\bibitem{KSSz97} J. Koml\'os, G.N. S\'ark\"ozy and E. Szemer\'edi (1997) 
Blow-up Lemma, Combinatorica, {\bf 17}, 109-123.

\bibitem{KSSz98} J. Koml\'os, G.N. S\'ark\"ozy and E. Szemer\'edi (1998), An
Algorithmic Version of the Blow-up Lemma, Random Struct. Alg., {\bf 12}, 297-312.


\bibitem{KS93} J. Koml\'os, M. Simonovits (1993), Szemer\'edi's Regularity Lemma and
its Applications in Graph Theory, Combinatorics, Paul Erd\H{o}s is eighty,
Vol. 2 (Keszthely, 1993), 295--352.

\bibitem{MM} Cs. Magyar, R. Martin (2006), Tripartite version of the Corr\'adi-Hajnal theorem, 
Discrete Mathematics, in preparation.

\bibitem{MSz} R. Martin, E. Szemer\'edi (2008), Quadripartite version of the Hajnal-Szemer\'edi theorem, to appear.

\bibitem{SZ76} E. Szemer\'edi (1976), Regular Partitions of Graphs, Colloques
Internationaux C.N.R.S N$^{\underline {\rm o}}$ 260 - Probl\`emes
Combinatoires et Th\'eorie des Graphes, Orsay, 399-401.

\end{thebibliography}
\end{document}